\documentclass[reqno,11pt]{article}
\usepackage{amsmath,amssymb}
\usepackage{amsthm,amscd,amsfonts,}
\usepackage{amsmath}
\usepackage{amsthm}
\usepackage{graphicx}
\usepackage{amsmath,amsthm}
\usepackage{amsmath,amstext,amsfonts}
\usepackage{subfigure}
\usepackage[T1]{fontenc}
\usepackage[utf8]{inputenc}
\usepackage{authblk}
\usepackage{breqn}

\newtheorem{theorem}{Theorem}[section]

\newtheorem{remark}{Remark}[section]

\newtheorem{definition}{Definition}[section]

\newcommand{\R}{\mathbb{R}}
\numberwithin{equation}{section}

\begin{document}

\title{\textbf{Complete Monotonicity and Inequalities of Functions Involving $\Gamma$-function}}
\author{M. Al-Jararha\thanks{mohammad.ja@yu.edu.jo}}
\affil{Department of Mathematics, Yarmouk University, Irbid, Jordan, 21163.}

\date{}
\maketitle

\vspace{0.2in}

\noindent \textbf{AMS (2000) Subject Classification}: 33B15.

\begin{abstract}
In this paper, we investigate the complete monotonicity of some functions involving gamma function. Using the monotonic properties of these functions, we derived some inequalities involving gamma and beta functions. Such inequalities arising in probability theory.
\end{abstract}
\maketitle
\section{Introduction} 
Completely monotonic functions play a major role in Probability and Mathematical Analysis due to their monotonic properties \cite{Alzer1,Alzer8,Miller}.
A function $f(z)$ is called completely monotonic on an interval $I\subset \R$ if it has derivatives of any order $f^{(n)}(z),\; n=0,1,2,3,\cdots$, and if 
\[
(-1)^nf^{(n)}(z)\geq 0 
\]
for all $ x\in I$ and all $n\geq 0.$  
Recently,  complete monotonicity  of functions involving gamma and beta functions have been considered in many articles, for example \cite{Guo,Qi,Giordano,Alzer3,Vogt,Ismail1}. Also, many articles have appeared to provide various inequalities of functions that involving  gamma  and beta functions, see \cite{Bustoz,Chen,Laforgia,Alzer2,Shabani,Sandor,Marcer,Garg,Agarwal,Alzer4,Alzer44,Alzer6,Dragomir,Alzer7,Nantomah,Maligranda,Batir,Ivady,Zhimin,Alzer9}. 
In this paper, we investigate the complete monotonicity of some functions involving gamma function. Using the monotonic properties of these functions, we derive some inequalities involving gamma function, where gamma function is defined by 
\[
\Gamma(x)=\int_0^\infty e^{-t}t^{x-1}dt, \; x>0.
\]
Due to the relation between gamma and beta function, we derived some inequalities that are involving beta function. Commonly,  beta function is defined by 
\[
\mathcal B(x,y)=\int_0^1t^{x-1}(1-t)^{y-1}dt,\; x\;\text{and}\;y>0,
\]
and its relation with the gamma function is given by 
\[
\mathcal B(x,y)=\displaystyle\frac{\Gamma(x+y)}{\Gamma(x)\Gamma(y)},\; x,y>0.
\]

The layout of the paper: In the first section, we prove our main results. In the second section, we apply  some of our main results to prove the inequality
$$\displaystyle\frac{\Gamma(x+1)\Gamma(x-a-b+1)}{\Gamma(x-a+1)\Gamma(x-b+1)}\displaystyle\frac{\Gamma(y-a+1)\Gamma(y-b+1)}{\Gamma(y+1)\Gamma(y-a-b+1)}\geq 1,\;y\geq x\geq a+b>b\geq a>0.$$
The last section is devoted for the concluding remarks.
\section{The Main Results}
In this section, we present and prove our main results. First, we present some useful definitions and theorems.  
\begin{definition}
A function $f(z)$ is called completely monotonic on an interval $I$ if it has derivatives of any order $f^{(n)}(z),\; n=0,1,2,3,\cdots$, and if 
\[
(-1)^nf^{(n)}(z)\geq 0 
\]
for all $ x\in I$ and all $n\geq 0.$ If the above inequality is strict for all $x\in I$ and all $n\geq 0,$ then $f(z)$ is called strictly completely monotonic.
\end{definition}

\begin{definition}
A function $f(z)$ is called logarithmically completely monotonic on an interval $I$ if its logarithm has derivatives $[\ln f(z)]^{(n)}$ of orders $n\geq 1$, and if 
\[
(-1)^n[\ln f(z)]^{(n)}\geq 0 
\]
for all $ x\in I$ and all $n\geq 1.$ If the above inequality is strict for all $x\in I$ and all $n\geq 1,$ then $f(z)$ is called strictly logarithmically completely monotonic.
\end{definition}

\begin{theorem}\label{thm1} \cite{Qi}. Every (strict) logarithmically completely monotonic function is (strict) completely monotonic. 
\end{theorem}
Now, we turn to prove our main results. 
\begin{theorem}\label{theorem1}
Let $a,b\geq 0$. Define $f(z)=\displaystyle\frac{\Gamma(z+1)\Gamma(z-a-b+1)} {\Gamma(z-a+1)\Gamma(z-b+1)},\; z>a+b-1.$ Then $f(z)$ is completely monotonic function.\label{thm2}
\end{theorem}
\noindent\emph{proof.} Let $f(z)=\displaystyle\frac{\Gamma(z+1)\Gamma(z-a-b+1)}{\Gamma(z-a+1)\Gamma(z-b+1)}$. Then $f(z)>0,$ and 
\begin{equation}\label{lneq}
\ln f(z)=\ln \Gamma(z+1)+\ln \Gamma(z-a-b+1) +\ln \Gamma(z-a+1)+\ln \Gamma(z-b+1).
\end{equation}
By differentiating Eq. \eqref{lneq}, we get
\begin{eqnarray}
\frac{d}{dz}\ln f(z)&=&\frac{\Gamma^\prime(z+1)}{\Gamma(z+1)}+\frac{\Gamma^\prime(z-a-b+1)}{\Gamma(z-a-b+1)}-\frac{\Gamma^\prime(z-a+1)}{\Gamma(z-a+1)}-\frac{\Gamma^\prime(z-b+1)}{\Gamma(z-b+1)}\nonumber\\
&=&\Psi(z+1)+\Psi(z-a-b+1)-\Psi(z-a+1)-\Psi(z-b+1),
\end{eqnarray}
where $\Psi(z)$ is the Digamma function (the logarithmic differentiation of the $\Gamma-$function). By using the following integral representation of Digamma function  $$\Psi(z)=-\gamma+\displaystyle\int_0^\infty\frac{(e^{-t}-e^{-zt})}{1-e^{-t}}dt,\; \; Re(z)>0,$$ where $\gamma =0.577218...$ is Euler's constant, we get
\begin{eqnarray}
\frac{d}{dz}\ln f(z)
&=&\Psi(z+1)+\Psi(z-a-b+1)-\Psi(z-a+1)-\Psi(z-b+1)\nonumber\\
&=&-\displaystyle\int_0^\infty\frac{e^{-t(z-a-b+1)}}{1-e^{-t}}\left(e^{-(a+b)t}-e^{-at}-e^{-bt}+1\right)dt\nonumber\\
&=&-\displaystyle\int_0^\infty\frac{e^{-t(z-a-b+1)}}{1-e^{-t}}\left(1-e^{-at}\right)\left(1-e^{-bt}\right)dt\nonumber\\
&\leq& 0,\;\;\; \forall z> a+b-1.
\end{eqnarray}
Therefore, $-\frac{d}{dz}\ln f(z)\geq 0, \;\forall z>a+b-1.$ Inductively, we have
\begin{equation}
(-1)^n\frac{d^n}{dz^n}\ln f(z)
=\displaystyle\int_0^\infty\frac{t^{n-1}e^{-t(z-a-b+1)}}{1-e^{-t}}\left(1-e^{-at}\right)\left(1-e^{-bt}\right)dt\geq 0,\;\; \forall z> a+b-1.
\end{equation}
Hence, by Theorem \ref{thm1}, $f(z)$ is completely monotonic function.
\begin{remark}\label{rem21}
Since $f(z)$ is completely monotonic. Then it is decreasing function. Hence, by using the asymptotic relation 
\begin{equation}\label{asymptoticrelation}
\displaystyle\lim_{z\rightarrow \infty}z^{(b-a)} \;\frac{\Gamma(z+a)}{\Gamma(z+b)}=1,
\end{equation}  
 we get , \begin{equation}\label{inq1}\displaystyle\frac{\Gamma(z+1)\Gamma(z-a-b+1)}{\Gamma(z-a+1)\Gamma(z-b+1)}\geq 1,\;\;z> a+b-1.
 \end{equation} 
 For a reference of the above asymptotic relation, see  eq. 13 in \cite{Mangus} (see also, \cite{Abramwtz,Tricomi}).
\end{remark}
\begin{theorem}\label{thm3}
Let $a_1,a_2,\cdots,a_n\geq0$, and let $\bar{a}=\sum_{i=1}^n a_i$. Then   $f(z)=\displaystyle\frac{\Gamma^{n-1}(z+1)\Gamma(z-\bar a+1)}{\prod_{i=1}^n\Gamma(z-a_i+1)},$ $\;z>\bar a-1$ is completely monotonic function.
\end{theorem}
\noindent\emph{proof.} Let $f(z)=\displaystyle\frac{\Gamma^{n-1}(z+1)\Gamma(z-\bar a+1)}{\prod_{i=1}^n\Gamma(z-a_i+1)},$ $\;z>\bar a-1$. Then
\begin{eqnarray}
(-1)^n\frac{d^nf(z)}{dz^n}&=&\displaystyle\int_0^\infty\frac{t^{(n-1)}e^{-t(z+1)}}{1-e^{-t}}\left((n-1)+e^{\bar at}-\sum_{i=1}^ne^{a_it}\right)dt\nonumber\\
&\geq&0,\;\;\;z>\bar a-1. \nonumber
\end{eqnarray}
This is correct since the function $\xi(t):=(n-1)+e^{\bar at}-\sum_{i=1}^ne^{a_it},\;t\geq 0,$ is nonnegative function. In fact, $\xi(0)=0$ and $\xi^\prime(t)=\bar ae^{\bar a t}-\sum_{i=1}^n a_ie^{a_it}=\sum_{i=1}^n a_ie^{\bar a t}-\sum_{i=1}^n a_ie^{a_it}=\sum_{i=1}^n a_i(e^{\bar a t}-e^{a_it})>0,$ since the exponential function is increasing function. Hence, $\xi^\prime(0)=0$ and $\xi^\prime(t)>0,\;\forall t>0.$ Therefore, $\xi(t)$ is nonnegative. Hence, by Theorem \ref{thm1}, $f(z)$ is completely monotonic.
\begin{remark}\label{rem22}
Let $a_1,a_2,\cdots,a_n\geq0$, and let $\bar{a}=\sum_{i=1}^n a_i$. Since   $$f(z)=\displaystyle\frac{\Gamma^{n-1}(z+1)\Gamma(z-\bar a+1)}{\prod_{i=1}^n\Gamma(z-a_i+1)},\;\forall z>\bar a-1$$ is completely monotonic function, and hence, it is decreasing function. Therefore, by using the asymptotic relation \eqref{asymptoticrelation}, 
we get
\[
f(z)=\displaystyle\frac{\Gamma^{n-1}(z+1)\Gamma(z-\bar a+1)}{\prod_{i=1}^n\Gamma(z-a_i+1)}\geq 1,\;\forall z>\bar a-1.
\]
\end{remark}
\begin{theorem}\label{thm4} Let $0\leq a_1\leq a_2\leq\cdots\leq a_n$, $0\leq b_1\leq b_2\leq \cdots\leq b_n$, and 
$m=\max\{a_1,a_2,\cdots,a_n,\\ b_1,b_2,\cdots,b_n\}.$ Moreover, assume that $\sum_{i=1}^k b_i\leq \sum_{i=1}^k a_i,\;k=1,2,\cdots,n$. Then $$f(z)=\displaystyle\frac{\prod_{i=1}^n\Gamma(z-a_i)}{\prod_{i=1}^n\Gamma(z-b_i)},\; z>m,$$ is completely monotonic function.
\end{theorem}
\noindent\emph{proof.} Let $f(z)=\displaystyle\frac{\prod_{i=1}^n\Gamma(z-a_i)}{\prod_{i=1}^n\Gamma(z-b_i)},\; z>m$. Then
\begin{eqnarray}
(-1)^n\frac{d^nf(z)}{dz^n}&=&\displaystyle\int_0^\infty\frac{t^{(n-1)}e^{-tz}}{1-e^{-t}}\left(\sum_{i=1}^ne^{a_it}-\sum_{i=1}^ne^{b_it}\right)dt\nonumber\\
&\geq&0,\;\;\;z> m. \nonumber
\end{eqnarray}
This inequality holds since $\sum_{i=1}^k b_i\leq \sum_{i=1}^k a_i,\;k=1,2,\cdots,n$ and the exponential function $e^x,\;x>0$ is convex and increasing function which implies that  $ \sum_{i=1}^ne^{a_it}\geq\sum_{i=1}^ne^{b_it}$ (see page 12 in \cite{maj}). Hence, by Theorem \ref{thm1}, $f(z)$ is completely monotonic function.

\begin{remark}\label{rem23}
 Let $0\leq a_1\leq a_2\leq\cdots\leq a_n$, $0\leq b_1\leq b_2\leq \cdots\leq b_n$, and 
$m=\max\{a_1,a_2,\cdots,a_n,\\ b_1,b_2,\cdots,b_n\}.$ Moreover, assume that $\sum_{i=1}^kb_i\leq \sum_{i=1}^k a_i,\;k=1,2,\cdots,n-1$ and $\sum_{i=1}^na_i=\sum_{i=1}^nb_i$. Then $f(z)$ is decreasing function as a result of the above theorem. Hence, by the limit \eqref{asymptoticrelation}, 
 we get  $\frac{\prod_{i=1}^n\Gamma(z-a_i)}{\prod_{i=1}^n\Gamma(z-b_i)}\geq 1,\; \forall z>m.$
\end{remark}

\begin{remark}\label{rem24} Let $0\leq a_1\leq a_2\cdots \leq a_n$ and let $b_i=\bar a:=\displaystyle\frac{\sum_{i=1}^na_i}{n},\; i=1,2,\cdots,n.$  Then
$\sum_{i=1}^kb_i=k\bar a=\frac{k}{n}\sum_{i=1}^k a_i<\sum_{i=1}^k a_i,\;k=1,2,\cdots,n-1$, and $\sum_{i=1}^na_i=\sum_{i=1}^nb_i$. Hence,  $f(z)=\displaystyle\frac{\prod_{i=1}^n\Gamma(z-a_i)}{\Gamma(z-\bar a)^n}$ is  completely monotonic $\forall z>\bar a $. Moreover, the inequality $\displaystyle\frac{\prod_{i=1}^n\Gamma(z-a_i)}{\Gamma(z-\bar a)^n}\geq 1$  holds $ \forall z>\bar a$. By letting $z=z+1$ in this inequality, we get that $\displaystyle\frac{\prod_{i=1}^n\Gamma(z-a_i+1)}{\Gamma(z-\bar a+1)^n}\geq 1,\; \forall z>\bar a-1.$ For a particular case if we let $n=2, \;a_1=x>0,\; a_2=y>,$ and $z=x+y$, then we get $\displaystyle\frac{\Gamma(x+1)\Gamma(y+1)}{\Gamma(\frac{x+y}{2}+1)^2}\geq 1,\; \forall x,y\geq 0.$ Particularly, we get $\displaystyle\frac{\Gamma(x)\Gamma(y)}{\Gamma(\frac{x+y}{2})^2}\geq \frac{(x+y)^2}{4xy},\; \forall x,y> 0.$ This inequality has been proved in different method in \cite{Ivady}.

\end{remark}


\begin{theorem}\label{thm5}
Let $a\geq 0$ and define $f(z)=\displaystyle\frac{\Gamma(z+a)\Gamma(z-a)}{\Gamma(z)^2},\;z>a.$ Then $f(z)$ is  completely monotonic function.
\end{theorem}
\noindent\emph{proof.} Let $f(z)=\displaystyle\frac{\Gamma(x-a)\Gamma(x+a)}{\Gamma(z)^2},\; z>a$. Then
\begin{eqnarray}
(-1)^n\frac{d^nf(z)}{dz^n}&=&\displaystyle\int_0^\infty\frac{t^{(n-1)}e^{-tz}}{1-e^{-t}}\left(e^{at}+e^{-at}-2\right)dt\nonumber\\
&=&\displaystyle\int_0^\infty\frac{t^{(n-1)}e^{-tz}}{1-e^{-t}}\left(e^{at}-1\right)\left(1-e^{-at}\right)dt\nonumber\\
&\geq&0,\;\;\;z> a. \nonumber
\end{eqnarray}
Hence, by Theorem \ref{thm1}, $f(z)$ is  completely monotonic function.
\begin{remark}\label{rem25}
Let $a\geq 0$. Then $f(z)=\displaystyle\frac{\Gamma(z+a)\Gamma(z-a)}{\Gamma(z)^2},\;z>a$ is decreasing function. By using the asymptotic relation  $\displaystyle\lim_{z\rightarrow \infty}z^{(b-a)}\frac{\Gamma(z+a)}{\Gamma(z+b)}=1,$ we get $\displaystyle\frac{\Gamma(z+a)\Gamma(z-a)}{\Gamma(z)^2}\geq 1,\;\forall z>a$. This inequality has been proved in \cite{Dragomir} by using classical integral inequalities. 
\end{remark}

\begin{remark}\label{rem26}
Using the same argument above, we can show that $$f(z)=\displaystyle\frac{\Gamma(z+a+1)\Gamma(z-a+1)}{\Gamma(z+1)^2},\;z>a-1$$ is completely monotonic, and so it is decreasing function. Consequently, the inequality 
\begin{equation}\displaystyle\frac{\Gamma(z+a+1)\Gamma(z-a+1)}{\Gamma(z+1)^2}\geq 1\label{inq2}\end{equation}
holds $\forall z>a-1$. Moreover, since $f(z)$ is decreasing function, then $f(z)\leq f(0)$. Hence, we get 
\begin{equation}\label{inq3}
1\leq \displaystyle\frac{\Gamma(z+a+1)\Gamma(z-a+1)}{\Gamma(z+1)^2}\leq \Gamma(1-a)\Gamma(1+a),\;0\leq a<1,\;z\geq 0.
\end{equation}
Moreover, by using the formula $\Gamma(z+1)=z\Gamma(z),$ we get
 $$\frac{z^2}{z^2-a^2}\leq \displaystyle\frac{\Gamma(z+a)\Gamma(z-a)}{\Gamma(z)^2}\leq \frac{z^2}{z^2-a^2}(\Gamma(1-a)\Gamma(1+a)),\;0\leq a<1,\;z>a.$$
 By using the formula $\Gamma(1-a)\Gamma(1+a)=\frac{\pi a}{\sin\pi a},\; 0<a<1$ (see page 48,\cite{Temme}), we get 
 $$\frac{z^2}{z^2-a^2}\leq \displaystyle\frac{\Gamma(z+a)\Gamma(z-a)}{\Gamma(z)^2}\leq \frac{\pi az^2}{\sin\pi a(z^2-a^2)},\;0<a<1,\;z>a.$$
  For the particular case $a=\frac{1}{2}$, we have  
 $$\frac{4z^2}{4z^2-1}\leq \displaystyle\frac{\Gamma(z+\frac{1}{2})\Gamma(z-\frac{1}{2})}{\Gamma(z)^2}\leq \frac{2\pi z^2}{4z^2-1},\;z>\frac{1}{2}.$$
 
Let $a=\frac{1}{2}$ in  \eqref{inq3}, then we have 
\[
1\leq \displaystyle\frac{\Gamma(z+\frac{3}{2})\Gamma(z+\frac{1}{2})}{\Gamma(z+1)^2}\leq \frac{\pi}{2},\;\;z\geq 0.
\]
Equivalently,
\[
\left(\frac{2}{2z+1}\right)^{\frac{1}{2}}\leq \displaystyle\frac{\Gamma(z+\frac{1}{2})}{\Gamma(z+1)}\leq \left(\frac{\pi}{2z+1}\right)^{\frac{1}{2}},\;\;z\geq 0.
\]
More precisely, we have
\[
\left(\frac{2z^2}{2z+1}\right)^{\frac{1}{2}}\leq \displaystyle\frac{\Gamma(z+\frac{1}{2})}{\Gamma(z)}\leq \left(\frac{\pi z^2}{2z+1}\right)^{\frac{1}{2}},\;\;z>0.
\]
\end{remark}

\section{Applications to Probability}
In this section, we prove the following inequality:  $$\displaystyle\frac{\Gamma(x+1)\Gamma(x-a-b+1)}{\Gamma(x-a+1)\Gamma(x-b+1)}\frac{\Gamma(y-a+1)\Gamma(y-b+1)}{\Gamma(y+1)\Gamma(y-a-b+1)}\geq 1, \;0<a\leq b<a+b\leq x\leq  y.$$ 
This inequality arises in the Probability in the problems dealing with the general binomial coefficients $\binom{\alpha}{\beta}=\frac{\Gamma(\alpha+1)}{\Gamma(\beta+1)\Gamma(\alpha-\beta+1)}.$ To prove this inequality, we present some useful remarks.

\begin{remark}\label{rem1}
Let $a,b> 0$  and let $M\geq a+b$. Then  $\forall z$ such that $0<a\leq b <a+b\leq z < M$, the function $\delta(z)=\displaystyle\frac{e^{-t(z-a-b+1)}}{1-e^{-t}}(1-e^{-at})(1-e^{-bt})$ is nonnegative and is not identically zero. Hence, as a consequence of  the proof of Theorem \ref{thm2}, we have  $f^\prime(z)<0.$  Therefore, $f(z)=\displaystyle\frac{\Gamma(z+1)\Gamma(z-a-b+1)} {\Gamma(z-a+1)\Gamma(z-b+1)}$ is strictly decreasing function on $ [a+b,M)$. Hence, $f(a+b)=\displaystyle\frac{\Gamma(a+b+1)}{\Gamma(a+1)\Gamma(b+1)}> 1.$ Let $g(z)=\displaystyle\frac{\Gamma(z-a+1)\Gamma(z-b+1)}{\Gamma(z+1)\Gamma(z-a-b+1)},\; z\geq a+b>b\geq a>0.$ Then $g(z)=\displaystyle\frac{1}{f(z)},$ and so $g(z)=\displaystyle\frac{\Gamma(z-a+1)\Gamma(z-b+1)}{\Gamma(z+1)\Gamma(z-a-b+1)}\leq 1,\; \forall z\geq a+b>b\geq a>0.$ Also, for $z\in [a+b,M)$, we have $g^\prime(z)>0.$

\end{remark}
\begin{remark}\label{rem3}
Let $0<a\leq b$. Define $h(x,y)=f(x)g(y),\;x,y\geq a+b,$ where $f(x)=\displaystyle\frac{\Gamma(x+1)\Gamma(x-a-b+1)}{\Gamma(x-a+1)\Gamma(x-b+1),}$ and $g(y)=\displaystyle\frac{\Gamma(y-a+1)\Gamma(y-b+1)}{\Gamma(y+1)\Gamma(y-a-b+1)}.$ Then $h(x,y)$ is positive and continuous on the rectangular domain $\mathcal D=[a+b,\infty)\times[a+b,\infty)$. Using the asymptotic relation $\displaystyle\lim_{z\rightarrow \infty}z^{(b-a)}\frac{\Gamma(z+a)}{\Gamma(z+b)}=1,$ we have $\lim_{\left\|(x,y)\right\|\rightarrow \infty}h(x,y)\rightarrow 1$.
\end{remark}

In the following theorem, we prove that $h(x,y)\geq 1$ on the triangular domain  $$\Omega=\left\{(x,y)\;|\; 0<a\leq b<a+b\leq x\leq y\right\} \subset\mathcal D.$$
\begin{theorem}
Let $h(x,y)=\displaystyle\frac{\Gamma(x+1)\Gamma(x-a-b+1)}{\Gamma(x-a+1)\Gamma(x-b+1)} \frac{\Gamma(y-a+1)\Gamma(y-b+1)}{\Gamma(y+1)\Gamma(y-a-b+1)}, \;(x,y)\in \Omega.$ Then $h(x,y)\geq 1, \; \forall (x,y) \in \Omega.$
\end{theorem} 

\noindent\emph{proof.} Since $h(x,y)$ is positive, continuous, and $\lim_{\left\|(x,y)\right\|\rightarrow \infty}h(x,y)\rightarrow 1$. Then, by using Remark \ref{rem1} and Remark \ref{rem3}, there exists a sufficiently large constant $M>0,$ such that $$h(a+b,M)=\displaystyle\frac{\Gamma(a+b+1)}{\Gamma(a+1)\Gamma(b+1)}\displaystyle\frac{\Gamma(M-a+1)\Gamma(M-b+1)}{\Gamma(M+1)\Gamma(M-a-b+1)}
\geq 1,$$
Now,  define the closed and bounded rectangular region:
 $$\mathcal R=\left\{(x,y)|\;0<a \leq b<a+b\leq x \leq y \leq M  \right\}.$$ Since $h(x,y)$ is continuous function on $\mathcal R$. Then it takes its absolute values on the boundaries of $\mathcal R$, or at the points in $\mathcal R$ where $\mathbf \nabla h(x,y)=\mathbf 0$ ($\nabla h(x,y)$ is the gradient of $h(x,y)$). Clearly, $\mathbf \nabla h(x,y)=f^\prime(x)g(y)\mathbf i+f(x)g^\prime(y)\mathbf j$, where $f(x)$ and $g(y)$ are defined in Remark \ref{rem3}. Moreover, by Remark \ref{rem1}, we have $f(x)>0,\;g(y)>0,\;g^\prime(y)>0,$ and $f^\prime(x)<0$ in $\mathcal R$. Hence, $\mathbf \nabla h(x,y)\neq \mathbf 0$ in $\mathcal R.$ Hence, the absolute values of $h(x,y)$ must occur at the boundaries of $\mathcal R.$ In fact, the boundaries of $\mathcal R$ are the line segments
\begin{enumerate}
\item $ l_1=\left\{(x,M)\;|\; a+b\leq x \leq M\right\},$
\item $ l_2=\left\{(a+b,y)\;|\; a+b\leq x \leq M\right\},$
and
\item $ l_3=\left\{(x,x)\;|\; a+b\leq x \leq M\right\}.$
\end{enumerate}
Obviously, $h(x,y)=h(x,x)=1$ on the line segment $ l_3.$ Moreover, on the line segment $l_1$, we have $$h(x,M)=\displaystyle\frac{\Gamma(x+1)\Gamma(x-a-b+1)}{\Gamma(x-a+1)\Gamma(x-b+1)}\displaystyle\frac{\Gamma(M-a+1)\Gamma(M-b+1)}{\Gamma(M+1)\Gamma(M-a-b+1)},\;a+b\leq x \leq M,$$ which is decreasing function in $x$ with a negative derivative. Similarly, on the line segment $l_1$, we have $$h(a+b,y)=\displaystyle\frac{\Gamma(a+b+1)}{\Gamma(a+1)\Gamma(b+1)}\displaystyle\frac{\Gamma(y-a+1)\Gamma(y-b+1)}{\Gamma(y+1)\Gamma(y-a-b+1)},\;a+b\leq y \leq M,$$ which is increasing function in $y$ with a positive derivative. Hence,  $h(x,y)$ takes its absolut values at the points $(a+b,a+b)$, $(M,M)$, and $(a+b,M)$. Clearly, $h(a+b,a+b)=h(M,M)=1.$ At the point $(a+b,M)$, we have
\begin{equation*}
h(a+b,M)=\displaystyle\frac{\Gamma(a+b+1)}{\Gamma(a+1)\Gamma(b+1)}\displaystyle\frac{\Gamma(M-a+1)\Gamma(M-b+1)}{\Gamma(M+1)\Gamma(M-a-b+1)}\geq 1.
\end{equation*}
This implies that $h(x,y)\geq 1, \; \forall (x,y) \in \Omega$. This completes the proof. 

\noindent Finally, we have the following remark:
\begin{remark}\label{rem33}
Let $\psi(x,y)=\displaystyle \frac{\Gamma(x+1)\Gamma(y-a+1)}{\Gamma(y+1)\Gamma(x-a+1)},\; y>x>a>0.$ Let $x\in(a,y).$ Then there exists $\alpha_x>0$, such that $x=y-\alpha_x.$ Hence, we define 
$$\xi(y)=\psi(y-\alpha_x,y)=\displaystyle \frac{\Gamma(y-\alpha_x+1)\Gamma(y-a+1)}{\Gamma(y+1)\Gamma(y-\alpha_x-a+1)},\; y>\alpha_x+a>a>0.$$ 
Then, by using the above Remark \ref{rem3}, we have $$\psi(y-\alpha_x,y)=\displaystyle \frac{\Gamma(y-\alpha_x+1)\Gamma(y-a+1)}{\Gamma(y+1)\Gamma(y-\alpha_x-a+1)}\leq 1,\;\; \forall y>\alpha_x+a>a>0.$$ Therefore, $$\psi(x,y)=\displaystyle \frac{\Gamma(x+1)\Gamma(y-a+1)}{\Gamma(y+1)\Gamma(x-a+1)}\leq 1,\;\; y>x>a>0.$$ Consequently, we also have $$\displaystyle \frac{\Gamma(x+1)\Gamma(y-a-b+1)}{\Gamma(y+1)\Gamma(x-a-b+1)}\leq 1,\; y>x>a+b>b\geq a>0.$$ 
\end{remark}
\section{Concluding Remarks}
We have seen in Remark \ref{rem21} that 
\begin{equation}\label{inq50}
\displaystyle\frac{\Gamma(z+1)\Gamma(z-a-b+1)} {\Gamma(z-a+1)\Gamma(z-b+1)}\geq 1,\; \forall z\geq a+b-1.
\end{equation}
Let $z=a+b,\; a,b\geq 0$ in \eqref{inq50}, then we get  
\begin{equation}\label{inq51}
\displaystyle\frac{\Gamma(a+b+1)}{\Gamma(a+1)\Gamma(b+1)}\geq 1, \; \forall a,b\geq0.
\end{equation} 
By using the fact that $\Gamma(z+1)=z\Gamma(z)$, we have
\begin{equation}\label{inq52}
\displaystyle\frac{\Gamma(a+b+1)}{\Gamma(a+1)\Gamma(b+1)}=\frac{(a+b)}{ab}\frac{\Gamma(a+b)}{\Gamma(a)\Gamma(b)}\geq 1,\;\forall a,b> 0.
\end{equation}
 Hence,
 \begin{equation}\label{inq53}
 \displaystyle\frac{\Gamma(a+b)}{\Gamma(a)\Gamma(b)}\geq \frac{ab}{(a+b)},\;\forall a,b>0.
 \end{equation}
 Since $\mathcal B(a,b)=\displaystyle\frac{\Gamma(a+b)}{\Gamma(a)\Gamma(b)}$, we get
 \begin{equation}\label{inq530}
 \mathcal B(a,b)\leq \displaystyle \frac{(a+b)}{ab},\; \forall a,b> 0
 \end{equation} 
 Moreover, by using the fact  $\displaystyle\mathcal B(x,y+1)=\frac{y}{x}\mathcal B(x+1,y)=\frac{y}{x+y}\mathcal B(x,y)$, we get
 \[
 \displaystyle\mathcal B(a,b+1)=\frac{b}{a}\mathcal B(a+1,b)=\frac{b}{a+b}\mathcal B(a,b)\leq \frac{1}{a},\;a,b> 0.
 \]
 Hence, $\displaystyle\mathcal B(a,b+1)\leq \frac{1}{a}$ and $\displaystyle\mathcal B(a+1,b)\leq \frac{1}{b}$ for $ a,b> 0.$
In \cite{Dragomir}, the authors  proved the following inequality:
\begin{equation} \label{inq531}
 \mathcal B(a,b)\leq \displaystyle \frac{1}{ab},\;  0<a,b\leq 1. 
\end{equation}
Clearly, if $0<a,b\leq 1$ and $a+b\leq 1$, the upper bound of $\mathcal b(a,b)$ given in \eqref{inq530} is better than the one given in \eqref{inq531}. Moreover, the inequality \eqref{inq530} is valid for $a,b>0$ while the inequality \eqref{inq531} is restricted on $a,b\in (0,1].$  

 Let $a=b$ in Eq. \eqref{inq1} and Eq. \eqref{inq3}.  Then we get
 \begin{equation}\label{inq54}
\displaystyle\frac{\Gamma(2a+1)}{\Gamma(a+1)^2}\geq 1, \;  a\geq 0,
\end{equation}
and 
\begin{equation}\label{inq55}
 \displaystyle\frac{\Gamma(2a)}{\Gamma(a)^2}\geq \frac{a}{2},\;a>0,
 \end{equation}
respectively. Assume that $0<a\leq 1.$ then we have the following inequality (see eq. 2.8 in \cite{Ivady}: 
\begin{equation}\label{inq550}
\displaystyle\frac{\Gamma(a)^2}{\Gamma(2a)}\geq \frac{2a-a^2}{a^2},
\end{equation}
By combining \eqref{inq55} with \eqref{inq550}, we get 
\[
\displaystyle\frac{2a-a^2}{a^2}\leq \frac{\Gamma(a)^2}{\Gamma(2a)}\leq \frac{2}{a},\; 0<a\leq 1.
\]    
%
 Assume that $a,b>0$ and $0<z\leq 1.$ Then inequality \eqref{inq50} implies that
 \begin{equation}\label{inq56}
 \displaystyle\frac{\Gamma(z-a-b+1)} {\Gamma(z-a+1)\Gamma(z-b+1)}\geq \frac{1}{\Gamma(z+1)}\geq 1, \; 0<a\leq b<a+b\leq z\leq 1.
 \end{equation}
 This is correct since $\Gamma(x+1)\leq 1, \forall z \in [0,1].$
 Moreover, we have
 \begin{equation}
 1\leq\displaystyle \frac{\Gamma(z+1)\Gamma(z-a-b+1)} {\Gamma(z-a+1)\Gamma(z-b+1)}\leq \displaystyle\frac{\Gamma(a+b+1)}{\Gamma(a+1)\Gamma(b+1)},\;z\geq a+b>b\geq a>0. 
 \end{equation}
Hence,  for $z\geq a+b>b\geq a>0$, we have 
 \begin{equation}\label{inq57}
 \displaystyle\frac{1}{\Gamma(z+1)}\leq \frac{\Gamma(z-a-b+1)} {\Gamma(z-a+1)\Gamma(z-b+1)}\leq \displaystyle\frac{1}{\Gamma(z+1)}\frac{\Gamma(a+b+1)}{\Gamma(a+1)\Gamma(b+1)}. 
 \end{equation}
 Therefore, $$\displaystyle \lim_{z\rightarrow \infty}\frac{\Gamma(z-a-b+1)} {\Gamma(z-a+1)\Gamma(z-b+1)}=0.$$ 
Now, let $a=b$ in \eqref{inq50}, then we get
 \[
 \displaystyle \frac{\Gamma(z+1)\Gamma(z-2a+1)} {\Gamma(z-a+1)^2}\geq 1, \; z>2a-1.
 \]
 Equivalently,
 \begin{equation}\label{inq58}
 \displaystyle \frac{\Gamma(z)\Gamma(z-2a)} {\Gamma(z-a)^2}\geq \frac{(z-a)^2}{z(z-2a)}=1+\frac{a^2}{z(z-2a)}, \; z>2a.
 \end{equation}
Clearly,
 \[
 \displaystyle \frac{\Gamma(z)\Gamma(z-2a)} {\Gamma(z-a)^2}\geq 1, \; z>2a.
 \]
In fact, in \cite{AnI} many lower bounds of $\frac{\Gamma(z)\Gamma(z-2a)} {\Gamma(z-a)^2}$ were presented. Some of these lower bounds are given in the following inequalities: 
\begin{equation}\label{inq580}
\frac{\Gamma(z)\Gamma(z-2a)} {\Gamma(z-a)^2}\geq \frac{a^2+z}{z},\;z>0\; \text{and}\; z-2a>0,
\end{equation}
and 
\begin{equation}\label{inq581}
\frac{\Gamma(z)\Gamma(z-2a)} {\Gamma(z-a)^2}\geq1+ \frac{a^2(z-2)}{(z-a-1)^2},\;z>2, z-2a>0,\; \text{and}\; a\neq 0,-1.
\end{equation} 



\end{document}